\documentclass[11pt]{article}
\usepackage{amsthm} 
\usepackage{amsmath} 
\usepackage{amssymb}
\usepackage{amsfonts} 
\usepackage{amscd}
\usepackage[mathscr]{eucal} 
\usepackage{verbatim}
\usepackage{graphics}

\newtheorem{thm}{Theorem}  \newtheorem*{mlem*}{Main Lemma}
\newtheorem*{thm*}{Theorem}  \newtheorem{prop}[thm]{Proposition}
\newcommand{\N}{{\mathbb N}}  \newcommand{\Z}{{\mathbb Z}}
 \newtheorem{lem}[thm]{Lemma}
\newtheorem*{cor*}{Corollary}

\begin{document} 
\title{Asymptotic formula for sum-free sets in abelian groups }
\author{R.Balasubramanian, \\ The Institute of
Mathematical Sciences\\ CIT Campus, Taramani\\ Chennai-600113,
India.\\ E-mail: balu@imsc.res.in\\
\\
Gyan Prakash\\
School of Mathematics, Harish Chandra Research Institute\\
Chatnag Road, Jhusi,\\
Allahabad - 211 019, India\\
E mail: gyan@mri.ernet.in  } \date{}
\maketitle 
\begin{abstract} 
Let $A$ be a subset of a finite abelian group $G$. We say that $A$ is sum-free
if there is no solution of the equation $x + y = z$, with $x, y, z$ belonging
to the set $A$. Let $SF(G)$ the sel of all sum-free subets of $G$ and 
$\sigma(G)$ deonotes the number $ n^{-1}(\log_2 |SF(G)|) $. 
In this article we shall improve the error term in asymptotic formula of
$\sigma(G)$ obtained 
in~\cite{GR}. The methods used are a slight refinement of methods of~\cite{GR}. 
\end{abstract}
Let $G$ be a finite abelian group of order $n$. A subset $A$ of $G$ is said to
 be sum-free if 
there is no solution of the equation $x + y = z$, with $x, y, z$ belonging
to the set $A$. Let $SF(G)$ denotes the set of all sum-free subsets of 
$G$. This article is motivated by the question of studying the cardinality
of the set $SF(G)$.\\
\\
\noindent
{\bf Definition:} \begin{enumerate}
\item Let $\mu(G)$ denotes the density of a largest sum-free
subset of $G$, so that any such subset has size $\mu(G)n$.\\
\item Given a set $B \subset G$ we say that $(x, y, z) \in B^3$
is a Schur
triple of the set $B$ if $ x + y = z$.
\end{enumerate}
Observing that all subsets of a sum-free set are sum-free we have the obvious 
inequality
\begin{equation} \label{obv}
|SF(G)| \geq 2^{\mu(G) n}
\end{equation} 
Let the symbol $\sigma(G)$ denotes the number $ n^{-1}(\log_2 |SF(G)|) $. Then 
from~\eqref{obv} it 
follows trivially that $\sigma(G) \geq \mu(G)$.\\
In this article we improve the results of Ben Green and Imre Ruzsa~\cite{GR}
and prove the following two results. 
The theorem~\ref{th} follows immediately from
theorem~\ref{BGschf} and a result from~\cite{GR}, namely
theorem~\ref{grane}. The methods used to prove theorem~\ref{BGschf}
are a slight refinements
of methods in~\cite{GR}.
\begin{thm} \label{th} Let $G$ be a finite abelian group of order $n$. Then
we have the following asymptotic formula
\begin{equation*}
\sigma (G) = \mu(G) + O(\frac{1}{(\ln n)^{1/27}} ).
\end{equation*}
\end{thm}
\begin{thm} \label{BGschf} 
There exist an absolute positive constant $\delta_0$ such that if $F \subset G$
as at-most $\delta n^2$ Schur triples, where $\delta \leq \delta_0$. Then 
\begin{equation} 
|F| \leq (\mu(G)   + C \delta^{1/3} )n
\end{equation}
where $C$ is an absolute positive constant.
\end{thm}
Earlier Ben Green and Ruzsa~\cite{GR} proved the following:
\begin{thm} (\cite{GR}, Theorem 1.8. ) \label{Gasym}
 Let $G$ be a finite abelian group of order $n$. Then
we have the following asymptotic formula
\begin{equation*}
\sigma (G) = \mu(G) + O(\frac{1}{(\ln n)^{1/45}} ).
\end{equation*}
\end{thm}
\begin{thm} \label{Gschf} (\cite{GR}, Proposition 2.2) Let $G$ be an abelian 
group, and 
suppose that $F \subseteq G$
has at-most $\delta n^2$ Schur triples. Then 
\begin{equation} 
|F| \leq (\mu(G)   + 2^{20} \delta^{1/5} )n
\end{equation}
\end{thm}
The following theorem is also proven in~\cite{GR}.
\begin{thm}(\cite{GR}, Corollary 4.3.) \label{12ml}
Let $G$ be an abelian group, and suppose that $F \subseteq G$
has at-most $\delta n^2$ Schur triples. Then 
\begin{equation}
|F| \leq (\text{max} (\frac{1}{3}, \mu(G)   + 3 \delta^{1/3} )n
\end{equation}
\end{thm}
The theorem~\ref{BGschf} follows immediately
from theorem~\ref{12ml} in the case 
$\mu(G) \geq \frac{1}{3}$. In the case $\mu(G) < \frac{1}{3}$, the 
theorem~\ref{BGschf} again follows immediately from theorem~\ref{12ml} in the 
case $\delta$ is not very ``small''. In the case $\delta$ is small we require
 Lemma~\ref{imp} where an estimate is done 
differently than in~\cite{GR}. For the rest of results we require to prove 
theorem~\ref{BGschf}, the methods used are completely identical as 
in~\cite{GR},
but the results used are not identical.\\
\\
\noindent 
For proving  theorem~\ref{th} we use the following result from~\cite{GR}.
\begin{thm}\label{grane} (\cite{GR}, Proposition 2.1') Let $G$ be an abelian 
group of cardinality $n$, 
where $n$ is sufficiently large. Then there is a family $\mathcal{F}$
of subsets of $G$ with the following properties
\begin{enumerate}
\item $\log_2|\mathcal{F}| \leq n(\ln n)^{-1/18}$;
\item Every $A \in $ SF$(G)$ is contained in some $F \in \mathcal{F}$;
\item If $F \in \mathcal{F}$ then $F$ has at-most $n^2(\ln n)^{-1/9}$
Schur triples.
\end{enumerate}
\end{thm}
\noindent
The theorem~\ref{th} follows immediately from theorem~\ref{grane} and
theorem~\ref{BGschf}. We shall reproduce the proof given in 
\cite{GR}. If $n$ is sufficiently large as required by 
theorem~\ref{grane} then associated to each $A \in SF(G)$ there is an
$F \in \mathcal{F}$ for which $A \subset F$. For a given $F$, the number of $A$
which can arise  in this way is at most $2^{|F|}$. Thus we have the bound
\begin{equation*}
|SF(G)| \leq \sum_{F \in \mathcal{F} } 2^{|F|} \leq | \mathcal{F}| 
max_{F \in \mathcal{F} } |F| 
\end{equation*}
Hence it follows that
\begin{equation}
\sigma(G) \leq \mu(G) + C\frac{1}{(\ln n)^{1/27}} + \frac{1}{\ln n^{1/18}}.
\end{equation}
But from the \eqref{obv} we have the inequality $\sigma(G) \geq \mu(G)$. 
Hence the 
theorem~\ref{th} follows.\\
\\
\noindent
In order to prove theorem~\ref{BGschf} we shall require 
the value of $\mu(G)$, which is now known for all finite abelian groups. 
In order to
explain the results we need the following definition.\\
\\
\noindent
{\bf Definition:} Suppose that $G$ is a finite abelian group of order
$n$. If  $n$ is divisible by any prime $p \equiv 2 {\rm (mod \ 3)}$
then we say that $G$  is type $I$. We say that G is type $I(p)$ if
it is type $I$ and if $p$ is the {\it least} prime factor of $n$ of the form
$3l + 2$. If $n$ is not divisible by any
prime $p \equiv 2 {\rm (mod \ 3)}$,  but $3|n$, then we say that $G$
is type $II$. Otherwise $G$ is said to be  type $III$. That is the group $G$
is said to be of
type $III$ if and only all the divisors of $n$ are congruent to $1$ modulo $3$.
\\
\\
The following theorem is due to P. H. Diananda and H. P. Yap~\cite{DY} for 
type~$I$ and type $II$ groups and due to  Green and Ruzsa~\cite{GR} for 
type~$III$ groups.
\begin{thm}(\cite{GR}, Theorem 1.5.) \label{mug}
Let $G$ be a finite abelian group of order $n$. Then the following holds.
\begin{enumerate}
\item If $G$ is of type $I(p)$ then $\mu(G) = \frac{1}{3} + \frac{1}{3p}$.
\item If $G$ is of type $II$ then $\mu(G) = \frac{1}{3}$.
\item If $G$ is of type $III$  then $\mu(G) = \frac{1}{3}  - \frac{1}{3m}$, 
where $m$ is the exponent of $G$.
\end{enumerate}
\end{thm}
\section{proof of theorem~\ref{BGschf} }
In case the group $G$ is not of type $III$ it follows from theorem~\ref{mug}
that $\mu(G) \geq \frac{1}{3}$ and hence the theorem~\ref{BGschf} follows 
immediate using theorem~\ref{12ml}. Therefore we are required to prove 
theorem~\ref{BGschf} for type $III$ groups only.\\
\\ 
For the rest of this article 
$G$ will be a finite abelian group of type $III$ and $m$ shall denote the 
exponent of $G$. The following proposition is
an immediate corollary of theorem~\ref{mug} and theorem~\ref{12ml}.
\begin{prop}\label{lm}
Let $G$ be an abelian group of type $III$. Let the order of $G$ be $n$ and the 
exponent of $G$ be $m$. If $F \subset G$
as at-most $\delta n^2$ Schur triples then 
\begin{enumerate}
\item $|F| \leq (\mu(G) + \frac{1}{3m} + 3\delta^{1/3})n$.
\item In the case $\delta^{1/3}m \geq 1$ then 
$|F| \leq (\mu(G) +  4\delta^{1/3})n$, that is the theorem~\ref{BGschf} holds
in this case.
\end{enumerate}
\end{prop}

Therefore to prove the theorem \ref{BGschf} we are left with the following
case.\\
{\bf Case:} The group $G$ is an abelian group of
type $III$,  order $n$ and exponent $m$. The subset $F \subset G$ has at most
$\delta n^2$ Schur triples and  $\delta^{1/3}m < 1$.\\
\\
Let  $\gamma$ be a character of $G$ and $q$ denotes the order of $\gamma$.
Given such $\gamma$  we define $H_j =
\gamma^{-1}(e^{2\pi ij/q}) $. We also denote the set $H_0 = ker(\gamma)$ by
just $H$. Notice that $H$ is a subgroup of $G$ and $H_j$ are cosets of $H$.
The cardinality of the coset $|H_j| = |H| = \frac{n}{q}$.
 The indices is to be
considered as residues modulo $q$ , reflecting the isomorphism $G/H
\equiv \Z/q\Z$. For any set $F \subset G$ we also define
 $F_j = F \cap H_j$ and $\alpha_j = |F_j|/|H_j|$.  

\begin{prop}Let $G$ be a finite abelian group of order $n$. Let $F$ be a
subset of $G$ having at~most $\delta n^2$ Schur triples where $\delta \geq 0$.
Let $\gamma$ be any character of $G$ and $q$ be its order. Also let $F_i$ and
$\alpha_i$ be as defined above. Then the following holds. \label{ineq}
\begin{enumerate}
\item If $x$ belongs to $F_i$ and $y$ belongs to $F_j$ then $ x + y$ belongs
  to $H_{i + j}$.
\item  The number of Schur triples $\{x, y, z\}$ of the set $A$
with $x$ belongs to $F_l, y$ belongs to $F_j$ and $z$ 
belongs to $F_{j + l}$ is at~least 
$|F_l|(|F_j| + |F_{j + l}| - |H|)$. In other words there are at~least 
$\alpha_l(\alpha_j + \alpha_{j + l} - 1) (\frac{n}{q})^2$ Schur triples  
$\{x, y, z\}$ of the set $F$
with $x$ belongs to the set $F_l$. 
\item Given any $l \in \Z/q\Z$ such  that $\alpha_i > 0$, it follows that 
for any $j \in \Z/q\Z$ the inequality 
\begin{equation}\label{alphal}
\alpha_j + \alpha_{j + l} \leq 
1 + \frac{\delta q^2}{\alpha_l}
\end{equation}
holds.
\item Given any $t \in \mathbb{R}$ we define  the set $L(t) \subset \Z/q\Z$ as 
follows. The set 
\begin{equation*}
 L(t) = \{i \in \Z/q\Z: \alpha_i + \alpha_{2i} \geq 1 +  t \}. 
\end{equation*}
Then it follows that
\begin{equation} \label{Lt}
\sum_{i \in L(t)} \alpha_i \leq  \frac{\delta q^2}{t}
\end{equation}
\end{enumerate}
\end{prop}

\begin{proof}
\begin{enumerate}
\item This follows immediately from the fact that $\gamma$ is an homomorphism.
\item In the case $|F_l|(|F_j| + |F_{j + l}| - |H|) \leq 0$, there is nothing
  to prove. Hence we can assume that the set $F_l \neq \phi$.
 Then for any $x$ which belongs to the set
$F_l$, the sets $x + F_j \subset H_{j + l}$. Since the set $F_{j + l}$ is also
a subset of $H_{j + l}$ and $|F_j| + |F_{j + l} |-|H| > 0$, it follows that
\begin{equation*}
|(x + F_j) \cap F_{j + l}| = ||F_j| + |F_{j + l} | - |(x + F_j) \cup F_{j + l}|
 \geq |F_j| + |F_{j + l} |-|H|.
\end{equation*}
Now for any $z$ belonging to the set $(x + F_j) \cap F_{j + l}$ there exist
$y$ belonging to $F_j$ such that $ x + y = z$. hence the claim follows.
\item From $II$ there are at~least 
$\alpha_l( \alpha_j + \alpha_{j + l} - 1)(\frac{n}{q})^2$ Schur triples of the
  set $F$. Hence the claim follows by the assumed upper bound on the number of 
Schur triples of the set $F$.
\item For any fixed $i \in L(t)$, taking $j = l = i$ in $II$, we get there are
  at least $\alpha_it(\frac{n}{q})^2$
 Schur triples $\{x, y , z\}$ of the set $F$ with
$x$ belonging to the set $F_i$. Now for given any two 
$i_1, i_2 \in L(t)$ such that 
$i_1 \neq i_2$ ,the sets $F_{i_1}$ and $F_{i_2}$ have no element in common. 
Therefore there are at least $t\sum_{i \in L(t)} \alpha_i$ Schur triples of
the set $F$. Hence the claim follows. 
\end{enumerate}
\end{proof}
Since the order of any character of an abelian group $G$ divides the order of 
group and $G$ is of type $III$, the order $q$ of any character $\gamma$ of
$G$ is odd and congruent to $1$ modulo $3$. Therefore $q = 6k + 1$ for some 
$k \in \N$. Let $I, H, M, T \subset \Z/q\Z$ denotes the image of natural
 projection
of the intervals
$\{k + 1, k + 2, \cdots, 5k - 1, 5k \},   \{k + 1, k + 2, \cdots, 2k - 1, 2k\},
 \{2k + 1, 2k + 2, \cdots, 4k - 1, 4k\}, \{4k + 1, 4k + 2, \cdots, 5k - 1, 
5k\} \subset \Z$ to $\Z/q\Z$. Then the set $I$ is divided into $2k$ disjoint
pairs of the form $(i, 2i)$ where $i$ belongs to the set $H \cup T$. 
  
\noindent
\begin{lem} Let $G$ be a finite abelian group of type $III$ and order $n$.
Suppose that $F \subset G$ has at-most $\delta n^2$ Schur triples. Let
$\gamma$ be a character of $G$. Let the order of $\gamma$ be equal to 
$q = 6k + 1$. 
Then the following inequality holds. \label{imp}
\begin{equation}
\sum_{i= k+1}^{5k} \alpha_i \leq 2k + 2\delta^{1/2}q^{3/2}
\end{equation}
\end{lem} 
\begin{proof}
The set $I = \{k+ 1, k + 2, \cdots, 5k\}$ is divided into $2k$ disjoint pairs
of the form $(i, 2i)$ where $i$ belongs to the set $H \cup T$.
Therefore it follows that 
\begin{equation}\label{IHT}
\sum_{i= k+1}^{5k} \alpha_i = \sum_{i \in H \cup T} (\alpha_i + \alpha_{2i} )
\end{equation}
Given a $t > 0$ we divide the set $H \cup T$ into two disjoint sets as
follows. We define the set
\begin{equation*}
 S = \{i \in H \cup T: \alpha_i + \alpha_{2i} \leq 1 + t \}
\end{equation*}
and
\begin{equation*}
L = \{i \in H \cup T: \alpha_i + \alpha_{2i} > 1 + t\}.
\end{equation*}
Therefore the sets $S$ and $L$ are disjoint and the set $H \cup T = S \cup L$.
Therefore it follows that 
\begin{equation}\label{SL}
 \sum_{i \in H \cup T} (\alpha_I + \alpha_{2i} ) =  \sum_{i \in S} (\alpha_i +
 \alpha_{2i} ) +  \sum_{i \in L} (\alpha_i + \alpha_{2i} )
\end{equation}
From proposition~\ref{alphal}~\eqref{Lt} we have the following inequality 
\begin{equation*}
\sum_{i \in L} \alpha_i \leq \frac{\delta q^2}{t}
\end{equation*}
Since for any $l \in \Z/q\Z$, the inequality $\alpha_l \leq 1$ holds trivially. 
It follows that
\begin{equation} \label{Li}
\sum_{i \in L} (\alpha_i + \alpha_{2i}) \leq  |L| + \frac{\delta q^2}{t}.
\end{equation}
Also the following inequality 
\begin{equation}\label{Si}
\sum_{i \in S} (\alpha_i + \alpha_{2i}) \leq  |S| + |S|t
\end{equation}
holds just by the definition of the set $S$.
Therefore from~\eqref{IHT},~\eqref{SL},~\eqref{Si},~\eqref{Li} it follows that
\begin{equation}
\sum_{i = k + 1}^{5k} \alpha_i \leq |L| + \frac{\delta q^2}{t} +  |S| + |S|t
\leq 2k + qt +  \frac{\delta q^2}{t}
\end{equation}
Now choosing $t = (\delta q)^{1/2}$ the lemma follows.
\end{proof}
{\bf Remark:} The sum appearing in last Lemma was estimated as 
$2k + \delta q^2$ in~\cite{GR}.
 There the estimate $\alpha_i + \alpha_{2i} \leq 
(\delta)^{1/2} q$ is used to estimate the 
right hand side of~\eqref{IHT}.\\ 
\\
Notice that Lemma~\ref{imp} holds for any character $\gamma$ of a group $G$
of type $III$. We would like to show that given $F \subset G$ having at~most 
$\delta n^2$ Schur triples and also assuming that
$(\delta)^{1/3} m < 1$ where $m$ is the exponent of $G$,
  there is a character $\gamma$ such that
$\alpha_l \leq C (\delta q)^{1/2}$ 
$i \in \{0, 1, 2, \cdots k\} \cup \{5k + 1, 5k + 2, \cdots, 6k\}$ where
$C$ is an absolute positive constant, $q$ is the
order of $\gamma$ and  $k = \frac{q - 1}{6}$. To be able to do this 
we recall the concept of special direction as defined in~\cite{GR}.
The method of proof of this part is completely identical as in~\cite{GR},
though the results are not.\\
\\
\noindent
Given any set $B \subset G$, and a character $\gamma$ of $G$ we define 
$\widehat{B}(\gamma) = \sum_{b \in B} \gamma(b)$. Given a set
$B \subset G$ fix a character $\gamma_s$ such that $Re\widehat{B}(\gamma)$ is 
minimal. We follow the terminology in~\cite{GR} and call $\gamma_s$ to be the
special direction of the set $B$. \\
The following Lemma is proven in~\cite{GR}, but we shall reproduce
the proof here for
the sake of completeness.
\begin{lem}(\cite{GR}, Lemma 7.1, Lemma 7.3. (iv)) \label{regam}
Let $G$ be an abelian group of type $III$. Let $F \subset G$ has at~most
$\delta n^2$ Schur triples.  Let $\gamma_s$ be a special direction
of the set $F$. Let $\alpha$ denotes the number $\frac{|F|}{|G|}$.
Then the following holds.
\begin{enumerate}
\item $Re\widehat{F}(\gamma_s) \leq \left(\frac{\delta}{\alpha (1 - \alpha)} 
-  \frac{\alpha^2}{\alpha (1 - \alpha)}\right)n. $
\item In case $\delta \leq \eta/5 $, then 
either $|F| \leq (\mu(G))n$ or the following inequality holds.
\begin{equation}
q^{-1} \sum_{j = 0}^{q - 1} \alpha_j \cos(\frac{2\pi j}{q}) + 
\frac{\mu(\Z/q\Z)^2}{1 -  \mu(\Z/q\Z)} < 6\delta.
\end{equation}
\end{enumerate}
\end{lem}

\begin{proof}
\begin{enumerate}
\item The number of Schur triples in the set $F$ is exactly 
$n^{-1} \sum_{\gamma}(\widehat{F}(\gamma))^2\widehat{F}(\gamma)$. This follows
  after the straightforward calculation, using the fact that
\begin{equation}\label{parseval}
\sum_{\gamma}\gamma(b) = 0 {\text if } b \neq 0, 
\end{equation}  
and is equal to $n$ if $ b = 0$ where $0$ here denotes the identity element
of
the group $G$. Therefore using the assumed
  upper bound on the number of Schur triples in the set $F$ it follows that
\begin{equation*}
n^{-1} \sum_{\gamma}(\widehat{F}(\gamma))^2\widehat{F}(\gamma) =
n^{-1} \sum_{\gamma \neq 1}(\widehat{F}(\gamma))^2\widehat{F}(\gamma)
+ n^{-1}\widehat{F}(1))^2\widehat{F}(1) \leq \delta n^2,
\end{equation*}
Where $\gamma = 1$ is the trivial character of the group $G$. Since 
$ n^{-1}(\widehat{F}(1))^2\widehat{F}(1) = (\alpha)^3 n^2$, it follows that
\begin{equation*}
Re\widehat{F}(\gamma_s) \sum_{\gamma \neq 1}(\widehat{F}(\gamma))^2 \leq
n^{-1} \sum_{\gamma \neq 1}(\widehat{F}(\gamma))^2\widehat{F}(\gamma)
\leq (\delta - \alpha^3) n^2
\end{equation*}
Since using~\eqref{parseval} it follows that
$\sum_{\gamma \neq 1}(\widehat{F}(\gamma))^2 = \alpha(1 - \alpha^2)n^2$, the
claim follows.
\item We have $Re \widehat{F}(\gamma_s) = 
|H|\sum_{j} \alpha_j \cos(\frac{2\pi j}{q})$.
Therefore in the case $|F| \geq \mu(G)$, then from (I) it follows that 
\begin{eqnarray}
q^{-1} \sum_{j = 0}^{q - 1} \alpha_j \cos(\frac{2\pi j}{q}) & \leq & 
\frac{\delta}{\alpha (1 - \alpha)} 
-  \frac{\alpha^2}{\alpha (1 - \alpha)}\\
q^{-1} \sum_{j = 0}^{q - 1} \alpha_j \cos(\frac{2\pi j}{q}) + 
\frac{\mu(G)^2}{1 -  \mu(G)} & \leq & 
\frac{\delta}{\alpha (1 - \alpha)}
\end{eqnarray}
Since from theorem~\ref{mug}  that $\mu(G) \geq  \mu(\Z/q\Z)$ it follows that
\begin{equation*}
\frac{\mu(G)^2}{1 -  \mu(G)} \geq \frac{(\mu(\Z/q\Z))^2}{1 -  \mu(\Z/q\Z)}.
\end{equation*}
The claim follows using this and the fact that $\mu(G) \geq \frac{1}{4}$,
which implies that $\frac{\delta}{\alpha (1 - \alpha)} \leq 6\delta$.
\end{enumerate}
\end{proof}
\begin{prop}\label{ls}
Let $G$ be an abelian group of type $III$. Let $n$ and $m$ denotes the order
and exponent of $G$ respectively.
Let $F \subset G$ has at~most $\delta n^2$ Schur triples 
 and $\delta^{1/3} m \leq 1$. Let $|F|\geq \mu(G)n$.
Let $\gamma_s$ be a special
direction of the set $F$ and $q$ be the order of $\gamma_s$. 
Let $ q = 6k + 1$ and $\alpha_i$ be as defined above. There exist an
positive absolute constants $q_0$ and $\delta_0$ such that if $q \geq q_0$
and $\delta \leq \delta_0$, then the following 
holds
\begin{equation}
\alpha_i \leq c (\delta q)^{1/2} \text { for all } i \in 
\{0, 1, \cdots, k - 1, k \} \cup \{5k + 1, 5k + 2, \cdots 6k -1\},
\end{equation}
where $c$ is an positive absolute  constant.
\end{prop}
\begin{proof}
If $F \subset G$ be the set as given, then $-F \subset G$ is also a set which
satisfies the same hypothesis as required in the statement of proposition.
It is also the case that $|F_j| = |(-F)_{-j}|$. Therefore to prove the 
proposition it is sufficient to show that
 \begin{equation*}
\alpha_i \leq c (\delta q)^{1/2} \text { for all } i \in 
\{0, 1, \cdots, k - 1, k \} 
\end{equation*}
for some positive absolute  constant $c$.\\
\noindent
Let $S = q^{-1} \sum_{j = 0}^{q - 1} \alpha_j \cos(\frac{2\pi j}{q}) + 
\frac{\mu(\Z/q\Z)^2}{1 -  \mu(\Z/q\Z)}$. Then from Lemma~\ref{regam} we have 
that 
\begin{equation} \label{cos}
S \leq 6\delta. 
\end{equation}
Now let for some $l \in \{0, 1, \cdots, k - 1, k\} $,  $\alpha_l > c (\delta
q)^{1/2}$ (where $c$ is a positive  number which we shall choose  later), then
we shall show that this violates~\eqref{cos}, provided $q$ and $c$ are
sufficiently large and $\delta$ is sufficiently small.  For this we shall find
the lower bound of  $M = q^{-1}\sum_{j = 0}^{q - 1} \alpha_j \cos(\frac{2\pi
j}{q})$.\\ Let $\gamma_j$ denotes  $\frac{(\alpha_j + \alpha_{j +
l})}{2}$. Then we have
\begin{equation*}
M = \frac{1}{q 2\cos(\frac{\pi l}{q}) } \sum_{j = 0}^{q - 1} \alpha_j\left
 (\cos\frac{(2j + l)\pi }{q} +  \cos\frac{(2j - l)\pi }{q}\right).
\end{equation*}
That is we have 
\begin{equation}
M = \frac{1}{q \cos(\frac{\pi l}{q}) }
 \sum_{j = 0}^{q - 1} \gamma_j
\cos\frac{(2j + l)\pi }{q} 
\end{equation}
Notice that $\cos(\frac{\pi l}{q})$ is not well defined if we consider
$l$ as an element of $\Z/q\Z$. This is because the function 
$\cos(\frac{\pi t}{q})$ as a function of $t$ is not periodic with period $q$ 
but is periodic with period $q$. But we have assumed that 
$l \in \{0, 1, \cdots, k - 1, k $, therefore 
the above computation is valid.\\
\\
\noindent
Since 
$ \delta^{1/2}q^{3/2} \leq  \delta^{1/2}m^{3/2} < 1$ is true by assumption, 
recalling Lemma~\ref{ineq} it follows that 
\begin{eqnarray}
2\gamma_j & = & \alpha_j + \alpha_{j + l} \leq 1 + \frac{1}{c}
\delta^{1/2}q^{3/2} \leq 1 +  \frac{1}{c}, \text{ for any } j \in  \Z/q\Z
\quad \\
\text{ and} \sum_{j} \gamma_j & = & \sum_{j} \alpha_j \geq \mu(G)n \geq 2k. 
\quad \label{lf}
\end{eqnarray}
\\
\noindent
The inequality~\eqref{lf} follows 
from the assumption that $|F| \geq \mu(G)n $.\\
\\
\noindent
Let $t_c $ denotes the number $1 + 1/c$. Let $E(c,q)$ denotes the minimum 
value of $\sum_{j = 0}^{q - 1} \gamma_j
\cos\frac{(2j + l)\pi }{q}$ subject to the constraints that 
$ 0 \leq \gamma_j \leq \frac{t_c}{2}$ and $\sum_{j} \gamma_j \geq 
2k
$.\\
The function $f: \Z \to \mathbb{R}$ given by 
$f(x) = \cos(\frac{(q + x)\pi }{q}) $ is an even function with period $2q$.
Also for $ 0 \leq x \leq q$ we have the following
\begin{equation}
f(0) < f(1) < f(2) < f(3) < \cdots < f(q - 1) < f(q)
\end{equation}
\noindent
Now to determine the minimum value of $E(c,q)$, we should choose $\gamma_j$
to be as large as we can when the function $\cos\frac{2j + l}{q}$ takes the 
small value.
 Now we have the two cases to 
discuss,
 the one when $l$ is even and  when $l$ is odd.
Now the image of function $g: \Z/qZ \to \mathbb{R}$ given by 
$g(j) =cos\frac{(2j + l)\pi }{q}$ is equal to 
$\{ f(x): x $  is even   \} 
in case 
$l$ is odd and is equal to  $\{ f(x): x$  is odd   \}
in case $l$ is even. From this it is also easy to observe that
the number of $j \in \Z/q\Z$ such that the function $\cos\frac{2j + l}{q}$
is negative is at~most $\frac{q + 1}{2}$.
Now let $ -\frac{q - 1}{2} - l \leq j \leq \frac{q - 1}{2} - l$ so that
$ -q \leq 2j + l \leq q$.
Now in case $l$ is odd then consider the case when $\gamma_j = \frac{t_c}{2}$
if
\begin{equation}
2j + l  = q - [\frac{k 
}{t_c} \frac{1}{2}] , \ldots, q -
2,q, q + 1, \ldots ,  q + [\frac{ k }{t_c} - \frac{1}{2}] \text{
and  } \gamma_j =  0 \text{ otherwise}.
\end{equation}
The condition $2[ \frac{k}{t_c} -1/2] + 1 \geq \frac{q + 1}{2}$ ensures that in
the above configuration for all possible negative values of
$ cos\frac{(2j + l)\pi }{q}$ the maximum possible weight $\frac{t_c}{2}$ is 
chosen.
This condition can be ensured  if $q \geq 11$ by choosing $c \geq c_1$ where
$c_1$ is sufficiently large positive absolute constant. 
Therefore after doing a small calculation one may check that for $c \geq c_1$
the following inequality 
\begin{equation} 
E(c,q) \geq - t_c 
\frac{\sin{\frac{2\pi[\frac{k}{t_c} -\frac{1}{2}]}{q}}}{2q\sin{\pi/q} 
\cos{\pi l/q}} - \frac{1}{q}
\end{equation}
holds.
In case $l$ is even and $c \geq c_1$ then choosing  $\gamma_j = \frac{t_c}{2}$
if
 \begin{equation}
 2j + l = q- [\frac{k}{t_c}] , \ldots, q - 1, q + 1, \ldots
,  q + [\frac{k}{t_c}] \text{ and weights } 0 \text{
otherwise},
\end{equation}  
we get that  the following inequality 
\begin{equation}
E(c) \geq - t_c \frac{\sin{\frac{2\pi[\frac{k}{t_c} ] + 1}{q}}}{2q\sin{\pi/q} 
\cos{\pi l/q}}
- \frac{t_c}{q}
\end{equation}
holds. Using this we get 
\begin{eqnarray}
S &\geq& -t_c \frac{\sin{\frac{2\pi[\frac{k}{t_c}]}{q}}}{2q\sin{\pi/q} 
\cos{\pi l/q}} + 
\frac{\mu(\Z/q\Z)^2}{1 -  \mu(\Z/q\Z)} \text{ when l is even and} \label{leven}\\
S &\geq& t_c 
\frac{\sin{\frac{2\pi[\frac{k}{t_c} -\frac{1}{2}]}{q}}}{2q\sin{\pi/q} 
\cos{\pi l/q}} - \frac{1}{q} + 
\frac{\mu(\Z/q\Z)^2}{1 -  \mu(\Z/q\Z)} \text{ when l is odd}. \label{lodd}\\
\end{eqnarray}
Now as $q \rightarrow \infty$ right hand side of~\eqref{leven} as well
as~\eqref{lodd} converges to the
\begin{equation*}
-t_c \frac{\sin\frac{2\pi}{3t_c}}{2\pi \cos\frac{\pi l}{q}} + 
\frac{1}{6}
\end{equation*}
Then let $\eta = 2^{-20}$, then choosing $c \geq c_2$ and $q \geq q_0$ we get 
and noticing that $l \leq \frac{q}{6}$ we get that
\begin{equation}
S \geq -\frac{1}{2\pi} + \frac{1}{6} - \eta = 8 \delta_0 \text{ say }. 
\end{equation}
The above quantity is strictly positive absolute constant. Then if
$\delta < \delta_0$, this contradicts~\eqref{cos}. Hence the Lemma follows. 
\end{proof}
To complete the proof of theorem~\ref{BGschf}, we require the following result 
from~\cite{GR}.
\begin{lem} (\cite{GR}, Proposition 7.2 )\label{sord}
Let $G$ be an abelian group of type $III$ and $n$ , $m$ be its order and 
exponent respectively. Let $F \subset G$ has at~most $\delta n^2$ Schur
triples, with $\delta^{1/3} m < 1$. Let $q$ be the order of special direction
such that $q \leq q_0$, where $q_0$ is a positive absolute constant as in 
Lemma~\ref{ls}. Also assume that $\delta \leq \frac{\eta}{q^5} = \delta'$, where 
$\eta= 2^{-50}$, then either $|F| \leq \mu(G)n$ or 
$\alpha_i \leq 64 \delta^{1/3} q^{2/3}$.
\end{lem}
Combining  Lemma~\ref{imp}, Lemma~\ref{ls} and lemma~\ref{sord} the 
theorem~\ref{BGschf} follows in the case $\delta^{1/3}m < 1$.
In the case $\delta^{1/3}m > 1$ the theorem follows from proposition~\ref{lm}. 

\end{document}